\documentclass[11pt,fleqn,reqno]{amsart}

\usepackage[english]{babel}
\usepackage{times,amsfonts,amsmath,amssymb,dsfont,mathrsfs,amsthm} 
\usepackage{pdfsync}
\usepackage{graphicx} 
\usepackage{mathabx}

\usepackage{color}

\definecolor{gr}{rgb}   {0.,   0.69,   0.23 }
\definecolor{bl}{rgb}   {0.,   0.5,   1. }
\definecolor{mg}{rgb}   {0.85,  0.,    0.85}
\definecolor{or}{rgb}   {0.9,  0.5,   0.}

\definecolor{webred}{rgb}{0.75,0,0}
\definecolor{webgreen}{rgb}{0,0.75,0}
\usepackage[citecolor=webgreen,colorlinks=true,linkcolor=webred]{hyperref}

\newtheorem{theorem}{Theorem}
\newtheorem{proposition}[theorem]{Proposition}

\theoremstyle{definition}

\theoremstyle{remark}

\newcommand{\beq}{\begin{equation}}
\newcommand{\eeq}{\end{equation}}

\newcommand{\R}{\mathbb{R}}
\newcommand{\C}{\mathbb{C}}

\newcommand{\rd}{{\mathrm d}}

\newcommand{\dS}{{\mathbb{S}}}

\newcommand{\ga}{\mathfrak{a}}

\newcommand{\gs}{\mathfrak{s}}
\newcommand{\gu}{\mathfrak{u}}

\newcommand{\gA}{\mathfrak{A}}

\newcommand{\gc}{\mathfrak{c}}

\newcommand{\bel}{\begin{equation} \label}
\newcommand{\ee}{\end{equation}}

\newcommand\eps{\varepsilon}

\renewcommand\div{\operatorname{div}}

\makeatletter
\def\section{\@startsection{section}{1}\z@{.9\linespacing\@plus\linespacing}%
  {.7\linespacing} {\fontsize{13}{15}\selectfont\scshape\centering}}
\def\paragraph{\@startsection{paragraph}{4}%
  \z@{0.3em}{-.5em}%
  {$\bullet$ \ \normalfont\itshape}}
\makeatother

\title[Extensions of the singular Robin Laplacian]{Self-adjoint and skew-symmetric extensions of the Laplacian with singular Robin boundary condition}

\author{Sergei A. Nazarov}
\address{Saint-Petersburg State University,
Universitetskaya nab., 7–9, St. Petersburg, 199034, Russia.}
\email{srgnazarov@yahoo.co.uk}

\author{Nicolas Popoff}
\address{Institut math\'ematique de Bordeaux, Universit\'e Bordeaux 1, UMR 5251, 33405 Talence CEDEX, France}
\email{nicolas.popoff@math.u-bordeaux1.fr}

\date{\today}

\begin{document}

\maketitle
\begin{abstract}
We study the Laplacian in a smooth bounded domain, with a varying Robin boundary condition singular at one point. The associated quadratic form is not semi-bounded from below, and the corresponding Laplacian is not self-adjoint, it has the residual spectrum covering the whole complex plane. We describe its self-adjoint extensions and exhibit a physically relevant skew-symmetric one. We approximate the boundary condition, giving rise to a family of self-adjoint operators, and we describe their eigenvalues by the method of matched asymptotic expansions. These eigenvalues acquire a strange behaviour when the small perturbation parameter $\varepsilon>0$ tends to zero, namely they become almost periodic in the logarithmic scale $|\ln \epsilon|$ and, in this way, ''wander'' along the real axis at a speed $O(\eps^{-1})$.
\end{abstract}
\subsection{Description of the singular problem}

In a domain $\Omega\subset \R^{2}$ enveloped by a smooth simple contour $\partial\Omega$, we consider the Laplacian with a Robin-type boundary condition $a\partial_{n}u-u=0$. Here, $a$ is a continuous function defined on $\partial\Omega$, and $\partial_{n}$ denotes the outward normal to $\partial\Omega$.

If $a$ is positive on $\partial\Omega$, the quadratic form $H^{1}(\Omega)\ni u\mapsto \|\nabla u;L^2(\Omega) \|^2-\|a^{-1/2}u,L^2(\partial\Omega) \|^2$ is naturally associated with this problem and, in view of the compact imbedding $H^{1}(\Omega) \subset L^{2}(\partial\Omega)$, this form is semi-bounded and closed, and thus defines a self-adjoint operator with compact resolvent. Therefore, the spectrum is an unbounded sequence of real eigenvalues accumulating at $+\infty$. Note that the first eigenvalue is negative, and goes to $-\infty$ if $a$ is a small positive constant, see \cite{LevPar08}.

Let $a$ become zero at a point $x_{0}\in \partial\Omega$. In this note, we will mainly consider the case where $a$ vanishes at order one, i.e. admits the Taylor formula 
 \bel{E:defa}
 a(s)=a_{0}s+O(s^2), \ \ s\to0 \quad \mbox{with} \ \ a_{0}>0,
 \ee
where $s$ is a curvilinear abscissa starting at $x_{0}$. For convenience, we denote $b_{0}:=a_{0}^{-1}$.

 Since we assume $a$ to be continuous, there should exist at least one other point where $a$ vanishes. However, our aim is to analyze the effect of one vanishing point of $a$, and the effect of several singular points can be explained along the same lines, thus we replace the problem by another one: we assume that $\partial\Omega$ is the union of two smooth curves $\Gamma_{1}$ and $\Gamma_{2}$ which meets perpendicularly, with $x_{0}$ in the interior of $\Gamma_{1}$, and that $a$ vanishes only at $x_{0}$, according to \eqref{E:defa}. We complete the Robin boundary condition on $\Gamma_{1}$ by a Neumann boundary condition on $\Gamma_{2}$. Therefore, our spectral problem is 
 \beq
\label{E:evproblem}
\left \{ 
\begin{aligned}
&-\Delta u =\lambda u \ \ \mbox{on}\ \  \Omega,
\\
&a\partial_{n}u-u=0 \ \ \mbox{on}\ \  \Gamma_{1}, \ \  \mbox{and} \ \  \partial_{n}u=0 \ \ \mbox{on}\ \  \Gamma_{2}.
\end{aligned}
\right. 
\ee
  The associated quadratic form is defined on 
  $D(q):=\{u\in H^{1}(\Omega), a^{-\frac{1}{2}}u_{|\Gamma_{1}}\in L^{2}(\Gamma_{1})\}$
   as follows
\begin{equation*}
D(q)\ni u\mapsto \int_{\Omega}|\nabla u|^2\rd x- \int_{\Gamma_{1}}a^{-1}|u|^2\rd s.
\end{equation*}
 It is not semi-bounded anymore. Thus, there is no canonical way for defining a self-adjoint operator associated with problem \eqref{E:evproblem}. The natural definition becomes the operator $A_{0}$ acting as $-\Delta$ on the domain
\bel{E:defA0}
D(A_{0}):=\{ u\in D(q), \Delta u\in L^{2}(\Omega), a \partial_{n}u-u=0\ \mbox{on}\ \ \Gamma_{1}, \partial_{n}u=0\ \mbox{on}\ \ \Gamma_{2}  \}.
\ee
Such a problem was studied in \cite{BerDen08,MarRoz09} in a model half-disk, for which the eigenvalue equation had the advantage to decouple in polar coordinates. The authors found that $A_{0}$ is non-self-adjoint. In \cite{MarRoz09}, they clarified the paradox from \cite{BerDen08} stating that, for any $\lambda\in \C$, problem \eqref{E:evproblem} has a nontrivial solution, by showing that the spectrum of $A_{0}^{*}$ is residual and coincides with the complex plane. 

\subsection{Goal and results}
In this note we explain how to find extensions of $A_{0}$ and give a better understanding of their spectrum, arguing with an asymptotical approach. We also exhibit a relevant skew-symmetric extension using a physical argument.



The domain of $A_{0}^{*}$ is 
\[D(A_{0}^{*}):=\{u\in L^{2}(\Omega), \Delta u \in L^{2}(\Omega),a \partial_{n}u-u=0\ \mbox{on}\ \Gamma_{1}, \, \partial_{n}u=0 \ \mbox{on}\ \Gamma_{2} \}.\]
To understand how far $D(A_{0}^{*})$ is from $D(A_{0})$, we exhibit two possible singular behaviors for functions in $D(A_{0}^{*})$ at the point $x_{0}$. Using Kondratiev theorem on asymptotics (\cite{Kon67}), we investigate a model problem in a half-plane and, as a result, describe $D(A_{0}^{*})$. We deduce, going over the domain $\Omega$, that the deficiency indices of $A_{0}$ are (1,1), and we classify its self-adjoint extensions 
using a parameter $\theta$ in the unit circle $\dS^{1}\subset \C$. The description of $A_{0}^{*}$ allows also to introduce a natural skew-symmetric extension of $A_{0}$ corresponding to a Sommerfield radiation condition at $x_{0}$.

Next, we approach our problem by a family of self-adjoint operators by choosing a suitable perturbation of the Robin coefficient $a$. This is done by means of the non-vanishing discontinuous function 
\bel{E:defaeps}
a_{\eps}(s)=a_{0}\mathrm{sign}(s)\eps+a(s)
\ee
 satisfying $\inf_{\Gamma_{1}}|a_{\eps}|=\eps$, and we study the discrete spectrum of the associated Robin Laplacian as $\eps\to0$. Using the method of matched asymptotic expansions, we find that its spectrum is related to eigenvalues of self-adjoint extensions, with a parameter $\theta_{\eps}$ oscillating in the logarithmic scale as $\eps\to0$. 

Finally, we describe the differences when the weight function satisfies $a(s)=a_{0}|s|+O(s^2)$ near the singular point, with $a_{0}>0$. In particular the number of singularities of functions in $D(A^{*}_{0})$ is now two or four, depending on whether $a_{0}>\tfrac{\pi}{2}$ or not. 

A similar result has been obtained in \cite{CheClaNaz16}, where an operator of the type $\div(\sigma \nabla)$ is considered in a bounded domain, where $\sigma$ is piecewise constant and changes sign along an interface crossing the boundary.

\subsection{Description of the adjoint operator}
\label{SS:modelproblem}
In this section, we investigate the following model problem in $\R^{2}_{+}$ : find $\gu\in L^{2}_{\mathrm{loc}}(\R^2_{+})$ such that
\beq
\label{E:evproblemHP}
\left \{ 
\begin{aligned}
&-\Delta \gu =0 \ \ \mbox{on}\ \  \R^{2}_{+},
\\
&\ga\partial_{n}\gu-\gu=0 \ \ \mbox{on}\ \  \partial \R^{2}_{+},
\end{aligned}
\right.
\ee
where $\ga$ is the first order approximation of $a$ near $x_{0}$: $\ga(s)=a_{0}s$.
Let $(r,\varphi)\in (0,+\infty)\times (-\frac{\pi}{2},\frac{\pi}{2})$ be the associated polar coordinates, the normal derivative reads $\partial_{n}\gu(s,0)=\mp r^{-1}\partial_{\varphi}\gu(r,\pm \frac{\pi}{2})$. 
The boundary condition is decoupled: The problem becomes in polar coordinates
\[
\left \{ 
\begin{aligned}
&-\partial_{r}^2\gu-r^{-1}\partial_{r}\gu-r^{-2}\partial_{\varphi}^2\gu =0 \ \ \mbox{on}\ \ (0,+\infty)\times (-\tfrac{\pi}{2},\tfrac{\pi}{2})
\\
&\forall r>0: \quad -a_{0}\partial_{\varphi}\gu-\gu=0 \ \ \mbox{at} \ \varphi=\pm \tfrac{\pi}{2}
\end{aligned}
\right.
\]
The spectrum of the transverse operator $-\partial_{\varphi}^2$ is given by solving the eigenvalue problem:
\beq
\label{E:evproblempolaireT}
-g''(\varphi)=\mu g(\varphi), \quad -a_{0}g'(\pm\tfrac{\pi}{2})-g(\pm\tfrac{\pi}{2})=0.
\ee
%
Its eigenvalues are $\{\mu_{k}, k\geq0\}:=\{-b_{0}^2\}\cup \{k^2, k=1,2,\ldots\}$. The eigenspace associated with $-b_{0}^2$ is generated by $g_{0}(\varphi)=e^{b_{0}\varphi}$, and the one associated with $k^{2}$ by $g_{k}(\varphi)=ka_{0}\cos(k\varphi)+\sin(k\varphi)$. We introduce two singular solutions of \eqref{E:evproblemHP}:
\beq
\label{E:D_gspm}
\gs^{\pm}(r,\varphi)=\gc r^{\pm i b_{0}}e^{b_{0}\varphi} \ \ \mbox{with} \ \ \gc=(1-e^{-2a_{0}^{-1}\pi})
\ee
where the choice for the normalizing factor $\gc$ will become clear in Proposition \ref{P:sympleticalgebra}.
 Note that $\gs^{\pm}\notin H^{1}_{\mathrm loc}(\overline{\R^{2}_{+}})$. 

  Let $\chi$ be a smooth cut-off function near $x_{0}$ with $\chi(x_{0})=1$ 
   and let $S^{\pm}$ be the functions deduced in $\Omega$ from $\gs^{\pm}$: $S^{\pm}(x)=\chi(x) \gs^{\pm}(r,\theta)$ through local polar coordinates $\Omega\ni x\mapsto (r,\theta) \in \R^{2}_{+}$ near $x_{0}$. 

As a consequence of the Kondratiev theorem on asymptotics (see \cite{Kon67} and, e.g. \cite[Ch.3]{NazPlam94}), we get

\begin{proposition}
\label{P:Kondratiev}
Let $u\in D(A_{0}^{*})$, then there exists $(c_{\mathrm{in}},c_{\mathrm{out}})\in \C^{2}$ such that
\bel{E:Kondratiev}
u=c_{\mathrm{in}}(u)S^{-}+c_{\mathrm{out}}(u)S^{+}+\tilde{u}
\ee
where 
$\tilde{u}\in H^{2}(\Omega)\cap D(q)$. 
Moreover, there exists $C>0$ such that, for all $u\in D(A_{0}^{*})$, we have 
\bel{E:estimeelliptique}
|c_{\mathrm{in}}(u)|+|c_{\mathrm{out}}(u)|+\|\tilde{u}\|_{L^2(\Omega)} \leq C(\|u\|_{L^2(\Omega)}+\|\Delta u\|_{L^2(\Omega)}).
\ee
\end{proposition}
This decomposition of the operator domain is sufficient to deduce the deficiency indices of the operator. On one hand, since the operator has real coefficient,
 $\dim(\ker(A^{*}_{0}+i))=\dim(\ker(A^{*}_{0}-i))$. On the other hand, standard decomposition together with last proposition implies $\dim(\ker(A^{*}_{0}+i))+\dim(\ker(A^{*}_{0}-i))=\dim(D(A_{0}^{*})/D(A_{0}))=2$. Therefore, $\dim(\ker(A^{*}_{0}\pm i))=1$, and the deficiency indices are (1,1). As a corollary, the spectrum of $A_{0}$ covers the whole complex plane.

\subsection{Self-adjoint extensions}
Once the domain of the adjoint is explicit, it is standard (\cite{Naz90,RofB69}) to find all self-adjoint extensions of $A_{0}$ by the use of the symplectic form 
\[\psi:(u,v)\mapsto \langle A_{0}^{*}u,v \rangle-\langle u,A_{0}^{*}v \rangle, \ \ \mbox{defined on} \ \ D(A_{0}^{*}).\]
As a consequence of integration by parts and symplectic algebra, we verify:
\begin{proposition}
\label{P:sympleticalgebra}
Let $u$ and $v$ in $D(A_{0}^{*})$, written in the form \eqref{E:Kondratiev}.
Then 
\[\psi(u,v)=i\left(c_{\mathrm{in}}(u)\overline{c_{\mathrm{in}}(v)}-c_{\mathrm{out}}(u)\overline{c_{\mathrm{out}}(v)}\right).\]
\end{proposition}

As a consequence of Proposition \ref{P:sympleticalgebra}, for any $u\in D(A_{0}^{*})$ , we obtain
\[\psi(u,u)=i(|c_{\mathrm{in}}|^2-|c_{\mathrm{out}}|^2).\]

Therefore, the self-adjoint extensions are the restrictions to the functions $u\in D(A_{0}^{*})$ for which $|c_{\mathrm{in}}(u)|=|c_{\mathrm{out}}(u)|$. We conclude with: 
\begin{theorem}
Let $\theta \in \dS^{1}$, and let $A_{0}(\theta)$ be the restriction of $A_{0}^{*}$ to the domain 
\[D(A_{0}(\theta))=\{u\in D(A_{0}^{*}),  c_{\mathrm{in}}(u)=e^{i\theta}c_{\mathrm{out}}(u) \}.\]
Then $A^{\mathrm{sa}}$ is self-adjoint extension of $A_{0}$ if and only if there exists $\theta\in \dS^{1}$ such that $A^{\mathrm{sa}}=A_{0}(\theta)$. 
\end{theorem}
From \eqref{E:estimeelliptique}, it follows that each domain of these extensions has continuous injection in $H^{2}$, which itself has compact injection in $L^{2}$. Therefore, each of these extensions has compact resolvent. Moreover, it is not semi-bounded from below, and we denote by $(\lambda_{k}(\theta))_{k\in \mathbb{Z}}$ the increasing sequence of eigenvalues of $A_{0}(\theta)$.

\subsection{The physical radiation condition and a skew-symmetric extension}

In link with the wave equation $-\partial_{t}^2W=-\Delta W$, analyzing propagation of the wave \[W^{\pm}(t,x):=e^{-i\sqrt{\lambda}\,t}S^{\pm}(x)=e^{-i\sqrt{\lambda}(t\mp\lambda^{-1/2}b_{0}\ln r)}g_0(\varphi),\] we can interpret the function $S^{-}$ as propagating from $x_{0}$, whereas $S^{+}$ would propagate toward $x_{0}$. 

For a fixed $\lambda\in \R$, the scattering theory (cf. \cite[Ch.5]{NazPlam94}) provides a solution of \eqref{E:evproblem} in the form 
\bel
{E:formsol}
\zeta_{\lambda}=S^{+}+e^{i\theta_{\lambda}}S^{-}+\tilde{\zeta}_{\lambda}, \quad \mbox{with} \ \ \tilde{\zeta}_{\lambda}\in D(q) \cap H^{2}(\Omega).
\ee
This solution is interpreted as the scattering wave associated with the entering wave $S^{+}$, and $e^{i\theta_{\lambda}}$ is the reflection coefficient, with $|e^{i\theta_{\lambda}}|=1$, according to conservation of energy.
%

Moreover, a natural skew-symmetric extension $\gA_{0}$ of $A_{0}$ can be defined in the domain
\[D(\gA_{0})=\{u\in D(A_{0}^{*}), c_{\mathrm{in}}(u)=0 \}.\]
This skew-symmetric extension corresponds to the natural Sommerfield radiation condition excluding entering waves.

\subsection{Wandering of the eigenvalues}

Assume that $\lambda_{k}(\theta)$ is a simple eigenvalue of $A_{0}(\theta)$, and denote by $C_{0}(S^{+}+e^{i\theta}S^{-})+\tilde{u}$ an associated normalized eigenfunction. Then standard computations show that $\lambda_{k}'(\theta)=-|C_{0}|^2$. Therefore an eigenvalue $\lambda_{k}(\cdot)$ is a non increasing function of $\theta\in \R$ wherever it is simple, moreover it is decreasing if $C_{0}\neq 0$. If $C_{0}=0$ for some $k$ and $\theta$, then the constant eigenvalue $\lambda_{k}(\theta)=\lambda$ is associated to a trapped mode, that is a non trivial solution of \eqref{E:evproblem} belonging to $D(q) \cap H^{2}(\Omega)$, which is in $D(A(\theta))$ for any $\theta\in \dS^{1}$.

The functions $\lambda_{k}(\cdot)$ are piecewise analytic, moreover they cannot be all constant, indeed in that case the range of all the eigenvalues $\lambda_{k}(\theta)$ would be a discrete set, which is a contradiction with the existence of the physical solution of \eqref{E:evproblem} of the form \eqref{E:formsol} for any $\lambda\in \R$.

Therefore, there exists at least one branch $\lambda_{k}(\cdot)$ which is decreasing where it is regular. This, combined with the periodicity of the spectrum, shows that 
\bel{E:unionR}
\bigcup_{(k,\theta)\in \mathbb{Z} \times \R} \lambda_{k}(\theta)=\R.
\ee


\subsection{The method of matched asymptotic expansions}

For $\eps>0$, we recall that $a_{\eps}$ was defined in \eqref{E:defaeps} as an approximation of $a$. Now the quadratic form
\[u\mapsto \int_{\Omega}|\nabla u|^2-\int_{\Gamma_{1}}a_{\eps}^{-1}(s)|u|^2\]
is well-defined and bounded from below in $H^{1}(\Omega)$. We denote by $A^{\epsilon}$ the corresponding self-adjoint operator. The strategy is now to construct quasi-modes for $A^{\epsilon}$ through eigenfunctions of an extension $A_{0}(\theta_{\epsilon})$, where $\theta_{\epsilon}$ is to be chosen. The result is 
\begin{theorem}
For all $k\geq0$ and for $\eps$ small enough, there exists $\theta_{\eps}\in \dS^{1}$ such that dist$(\lambda_{k}(\theta_{\eps}),\sigma(A^{\eps})) \to 0$ as $\eps\to0$. Moreover, the mapping $\eps\mapsto \theta_{\eps}$ is periodic with respect to $\ln \eps$, and 
\bel{E:wanderingev}\R = \bigcup_{\eps\downarrow0}\sigma(A^{\eps}).
\ee
\end{theorem}
The procedure is as follows: 
\paragraph{Far-field expansion:} Outside a fixed neighborhood of $x_{0}$, take a function $u^{out}$ as an eigenfunction of $A_{0}(\theta_{\epsilon})$, where $\theta_{\eps}$ is to be chosen. Therefore it behaves near $x_{0}$ as follows: 
\[u^{out}:(r,\varphi)\mapsto C( r^{ib_{0}}+e^{i\theta_{\eps}} r^{-ib_{0}})e^{a_{0}^{-1}\varphi}+\tilde{u}^{out},\]
where $\tilde{u}^{out}$ is regular and small near 0.

\paragraph{Near-field expansion} In local coordinates near $x_{0}$, we perform the scaling $x=\eps \xi$, and considering bounded eigenvalues, we get to solve \eqref{E:evproblemHP} with $\ga(\xi_{1}):=a_{0}(\mathrm{sign}(\xi_{1})+\xi_{1})$. In order to investigate the behavior of solutions to this problem at infinity, we perform the inversion $\xi\mapsto\eta=|\xi|^{-2}\xi$, which leads to the behavior at 0 of \eqref{E:evproblemHP} but with the weight function $a_{0}:\eta_{1}\mapsto\eta_{1}+\mathrm{sign}(\eta_{1})\eta_{1}^2$ in the boundary condition. 
Near the origin $\eta=0$, we can neglect the part $\mathrm{sign}(\eta_{1})\eta_{1}^2$, and according to the Kondratiev theory, there exists a solution of such a problem which behaves as 
\[ \eta \mapsto  \gs^{-}(\eta)+e^{i\theta}\gs^{+}(\eta) +O(|\eta|), \ \ \theta\in \dS^{1} \ \mbox{fixed}.\]
Therefore, we obtain a solution of \eqref{E:evproblemHP} with weight $\ga$, which produces after rescaling a solution of the Laplace equation in $\Omega$ which behaves in a neighborhood of $x_{0}$, as 
\[u^{in}:(r,\varphi)\mapsto \tilde{C}(\eps^{ib_{0}} r^{-ib_{0}}+e^{i\theta}\eps^{-ib_{0}} r^{ib_{0}})e^{a_{0}^{-1}\varphi}+\tilde{u}^{in}\]
where $\tilde{u}^{in}$ is decaying outside the neighborhood, and $\tilde{C}$ is a normalization factor. 

\paragraph{Matching expansions and conclusion} Matching the two previous extensions we obtain 
\bel{E:thetaeps}
\theta_{\eps} = \theta - 2b_{0}\ln \eps\pmod{2\pi}.
\ee
This formal approach is validated by constructing the quasi-mode from the previous ans\"atze using cut-off functions: define
\[u^{as}=\chi^{in}u^{in}+\chi^{out}u^{out}-\chi^{in}\chi^{out}C(\gs^{+}+e^{i\theta_{\eps}}\gs^{-}),\]
 where $\chi^{in}$ (respectively, $\chi^{out}$) is localized in a bounded neighborhood of $x_{0}$ (respectively, outside a neighborhood of $x_{0}$ of size $\epsilon$). Evaluating $(A^{\eps}-\lambda_{k}(\theta_{\eps}))u^{as}$, we get that $\lambda_{k}(\theta_{\eps})$ is close to the spectrum of $A^{\eps}$ for $\eps$ small enough. Note that $\theta_{\eps}$ is periodic w.r.t. $\log \eps$ and runs over $\dS^{1}$ at the rate $O(\eps^{-1})$ as $\eps\to0$. Then, \eqref{E:wanderingev} follows from \eqref{E:unionR}.

\subsection{Further questions}
When the weight function satisfies $a(s)=a_{0}|s|+O(s^2)$, with $a_{0}>0$, the situation depends on the parameter $a_{0}$ as described here: The transverse operator in the angular variable in the model half-plane $\R^{^2}_{+}$ is still $-\partial_{\varphi}^2$, but the boundary condition at $\varphi=\tfrac{\pi}{2}$ in \eqref{E:evproblempolaireT} now becomes $a_{0}g'(\tfrac{\pi}{2})-g(\tfrac{\pi}{2})=0$. The negative spectrum of this operator depends on $a_{0}$ as follows: 
\begin{enumerate}
\item If $a_{0}>\frac{\pi}{2}$, then there is one negative eigenvalue, and the others are positive. It produces two oscillatory solutions. 
\item If $a_{0}=\frac{\pi}{2}$, there is one negative eigenvalue, and null is also an eigenvalue. There are two additional solutions for the model problem, one has the form $g_{0}(\varphi)\ln r$ and the other is contant w.r.t. $r$.
\item If $a_{0}<\frac{\pi}{2}$, then there are two negative eigenvalues. They produce four oscillatory solutions. 
\end{enumerate}
Situation (1) can be analyzed exactly in the same way we described here. However situations (2) and (3) are much more different. In particular, the deficiency indices are (2,2), and self-adjoint extensions are parameterized by two-by-two unitary matrices. The method of matched asymptotic expansions does not provide a parameter extension explicit as in \eqref{E:thetaeps}, but a family of unitary matrices depending on $\eps$ (\cite{Naz99}). This family does not always coincide with the set of all unity matrices as $\eps\to0$, but it is sufficient for the construction of approximations.

\subsection*{Acknowledgement}
This work is performed within grant 17-11-01003 of Russian Science Fondation.

\bibliographystyle{abbrv}
\bibliography{bibliopof}

\end{document}